\documentclass{article}

\def\s{\rule{0in}{.35in}}

\def\var{\mbox{\rm var}}
\def\cov{\mbox{\rm cov}}

\usepackage{graphicx}
\usepackage{pslatex}

\usepackage{url}

\begin{document}

\title{On Steel's Test with Ties}
{\author{F.W. Scholz\thanks{Email: fscholz@u.washington.edu}}

\maketitle
\begin{abstract}
This note revisits Steel's multiple comparison test which uses Wilcoxon statistics in pairwise 
comparisons of several treatment samples with a common control sample. It
derives means, variances and covariances of the Wilcoxon statistics under the 
conditional randomization distribution, given the tie pattern in the pooled samples. 
Sample sizes do not have to be equal. 
Under the randomization distribution asymptotic multivariate normality of 
these Wilcoxon statistics is established.
This widens the scope of normal approximations to conditional tests,
assuming independent samples of respective
sizes $n_0, n_1,\ldots,n_K$ from any common population or randomized treatment assignment
to $N=n_0 + n_1+\ldots+n_K$ experimental subjects.
Significance probabilities are obtained using a normal approximation and a single quadrature.
In the continuous shift model the simultaneous tests are 
converted to simultaneous confidence bounds and intervals.
This is implemented in the {\sf R} package {\tt kSamples}. 
Extensions to all pairwise Wilcoxon test comparisons are 
discussed.

\begin{flushleft}
{\bf Keywords}: Nonparametric Methods, Multiple Comparisons, Multivariate Analysis, Normal Approximation

{\bf AMS Subject Classifications}: 62G10, 62G99
\end{flushleft}
\end{abstract}

\section{Introduction and Review}
This note revisits Steel's multiple comparison test which uses Wilcoxon test statistics in pairwise 
comparisons of several treatment samples with a common control sample. It is usually assumed that 
the treatment sample sizes are the same but possibly different from the control sample size.
It was also assumed by \cite{Steel59} that the sampled distributions are continuous 
in order to avoid ties. He provided tables for a limited set of
scenarios under these assumption. The main reason for avoiding ties was the intractability of the null distribution, since it would have to be dealt with conditionally
given the tie pattern which quickly become numerous.
\cite{Steel60} makes some suggestions concerning the handling of ties.
When they occur in the same sample, Steel suggests to ignore them, i.e., rank them in any consecutive order. 
This runs counter to the traditional treatment of ties which uses the conditional 
randomization distribution
given the tie pattern.
This was pointed out by \cite{Lehmann06} (p. 19) in discussing the treatment
of ties in the Wilcoxon test.
When ties spread over different samples, Steel suggests to rank them in consecutive order
to achieve a maximum or minimum value of the test criterion, whichever way would 
lead to a conservative procedure. This allows the use of tables or algorithms
based on the assumption of no ties
and could be viewed as a conservative 
tie breaking procedure, assuming that ties are due to rounding.

Algorithms for the exact
calculations of significance probabilities and close bounding values for the same were given in \cite{vdWiel02}. These are also applicable in the presence of ties.
However, computation time can still be forbidding when sample sizes or the number of samples are not small. Also, it is 
assumed that the treatment samples all have the same size.
The normal approximation comparisons in \cite{vdWiel02}
make use of the
variance and covariance formulas given in \cite{Miller81}, which assume no ties in the data.

As pointed out in \cite{Miller81}, a large sample approximation for the null distribution
of Steel's test statistic in the continuous population sampling case is implicitly  given by \cite{Lehmann63} and in a much more general context by \cite{MunzelHothorn01},
where sampling also from discrete distributions and under non-null situations is considered.  
However, when sampling a common discrete distribution (under control and any of
the treatments) the nonparametric nature of the test is lost, since the nature of
the common discrete sampled population has a strong effect on the null distribution, 
as pointed out in \cite{Lehmann06}, Chapter 1, Section 4. Only asymptotically
the nonparametric nature of the test can be recovered by estimating the 
covariance matrix of the involved statistics.  This is similar to the null distribution
of the ordinary $t$-test 
being asymptotically 
nonparametric when the variances of the sampled distributions are finite.

In the presence of ties one can salvage the nonparametric nature of rank tests by using midranks
of the pooled data and viewing
the distribution of the midranks conditionally given the tie pattern. This pattern is captured by the 
$N$ pooled ordered sample values $Z_{(1)}\leq Z_{(2)}\leq \ldots \leq Z_{(N)}$, 
or equivalently by the $e$
distinct ordered pooled observations $u_1< u_2 < \ldots < u_e$ together with their observed multiplicities
$d_1, \ldots, d_e$ as caused by the ties.
This conditional distribution is induced by the assumed random sampling
from a common distribution. It consists of all possible ways of splitting the pooled sample values
into $K+1$ samples of respective sizes $n_0, n_1, \ldots, n_K$, with 
$N=n_0+n_1+\ldots+n_K$, each one of these splits
having probability ${N\choose n_0,n_1,\ldots, n_K}^{-1}$. This distribution 
also arises naturally in the context of dealing with $N$ experimental units
from which observations are to be gathered. These units are randomly divided 
into $K+1$ groups of sizes $n_0, n_1, \ldots, n_K$, each group being administered 
a control (placebo) or one of the $K$ treatments. If there is no difference in effects between
the control and any of the treatments, then we have the same probability
distribution as before, i.e., probability ${N\choose n_0,n_1,\ldots, n_K}^{-1}$ for each split of the 
midranks into $K+1$ groups. However, any conclusions that are drawn from such an analysis
pertain to just the $N$ involved experimental units. Whether such conclusions can be generalized
based on the assumption that the $N$ units could be considered a random sample from a 
larger population is usually mostly a step of faith.

The monograph by \cite{Edgington07} makes a very strong case for the use of
randomization tests in contrast to the often invoked assumption of dealing with 
true random samples from some posited populations. 
The random sample assumption is often not satisfied in practice
while randomization with any given set of experimental subjects (aside from ethical concerns) is completely 
under the control of the experimenter and is quite manageable. It provides
a valid framework for any probability calculations, such as significance probabilities
or $p$-values. It is also easily 
grasped by the layperson.
This  distinction
is also emphasized in the treatment of rank tests in \cite{Lehmann06} by giving separate discussions
to these two approaches, introducing the randomization approach first for most cases covered.
This two-track approach may have been partly for didactic reasons in trying to reach a wider audience.

In the case of Steel's test this approach was not followed in \cite{Lehmann06}
for several reasons. Exact tables
were quite limited, see \cite{Steel59}. At the original publication of
\cite{Lehmann75} the large sample approximation was available only for sampling from 
continuous populations which also serves
for the randomization approach when there are no ties, see \cite{Miller81}.
When there are ties, no limit theory for the randomization approach existed.
\cite{MunzelHothorn01} provide such limit theory for the population random sampling situation
but it does not carry over to the randomization model, which conditions on the tie pattern
in the pooled samples. 

Finally, in 1975
software for estimating  $p$-values based the conditional
randomization distribution via simulation
was not yet widely available. This all has changed since then.
With the addition of the Steel test to the \emph{R} package
{\tt kSamples}, see \cite{Scholz16}, simulation is easily done.
Also provided there is a normal approximation in the presence of ties
under any configuration of sample sizes. 
The purpose of this note it to provide the theoretical 
underpinnings used in that package. For a very limited set of sample size
scenarios {\tt kSamples} also permits exact calculation of $p$-values in the presence of ties.

\section{Assumptions and Notation}
Suppose we are dealing with $n_0$ control observations
and  $n_i$ observations under treatment $i$ for $i=1,\ldots,K$.  It is assumed that the treatments and controls
are assigned by randomization to the $N=n_0+n_1+\ldots+n_K$ subjects. 
Each such partition into $K+1$ groups
has the same probability 
${N\choose n_0,n_1,\ldots, n_K}^{-1}$ under the null hypothesis
$H$ of no treatment effects.
The null distribution of any statistic
based on these $K+1$ samples with resulting observations $w_1,\ldots, w_N$
is a direct consequence 
of this equal probability randomization model.

Alternatively, we may assume that we deal with
$K+1$ independent random samples from $K+1$ populations.
Under the null hypothesis $H^\prime$ that the sampled populations are identical
we can view these samples as having arisen from an equivalent two stage
process. First draw one sample of size $N$ from that common population
and then randomly split that sample into $K+1$ subsamples of respective sizes
$n_0, n_1, \ldots, n_K$,
designating the various groups as 
respective control and treatment samples. The second stage is equivalent to our previous 
randomization model, conditionally given the pooled sample $w_1,\ldots, w_N$.

Since we will be dealing only with 
rank statistics we assume that these observation consist just of their midranks.
Assume that the midranks $w_1,\ldots, w_N$,
consists of $e$ distinct values
$u_1 < u_2 < \ldots < u_e$ with respective multiplicities $d_1,\ldots, d_e$,
with $d_1+\ldots+d_e=N$. By ${\cal T}=\{e,u_1<\ldots<u_e, d_1,\ldots,d_e\}$ we indicate the 
tie pattern in the observed midrank vector $w=(w_1,\ldots,w_N)$.

We compute the Wilcoxon rank sum test statistics
$W_{i}^\ast$ for each sample pair of $n_0$ and $n_i$ observations, $i=1,\ldots,K$,
i.e., $W_{i}^\ast$ is the sum of midranks of the $i^{\rm th}$ treatment sample
when re-ranking the $n_0+n_i$ observations in that combined group.
The asterisk in $W_i^\ast$ indicates the extension to tied 
observations as in \cite{Lehmann06}.
Under the 
hypothesis $H$ of no treatment effect, i.e., all $K+1$ sampled populations are the same or the randomly assigned treatments have no effect, we are dealing with a conditional null
distribution of $K$ correlated
Wilcoxon statistics. The correlation is induced by the common control sample used in the respective comparisons and conditioning is on
${\cal T}$, using as null distribution the permutation distribution
induced either by the random assignments of treatments and control
or by population sampling in conjunction with conditioning on the observed tie pattern.

Assuming $n_1=\ldots=n_K$, \cite{Steel59}  proposed 
\[
S = \max(W_1^\ast,\ldots,W_K^\ast)
\]
to test $H$ against the alternative that 
at least some of the treatments have a positive effect so that large values
of $S$ would be grounds for rejecting $H$ in favor of that alternative. 
Aside from making the special assumption of equal treatment sample sizes
a continuous population sampling model was also assumed. This avoids the problem of ties.
In the case of tied 
observations, as they are bound
to occur in practice, we will use midranks based on the pooled sample values. 

When sample sizes across treatments are not the same, either by design or because of completely at random missing data, 
a reasonable extension of the 
above test statistic $S$ is
\[
\hat{S}_{\max} =\max\left(
\frac{W_1^\ast-\mu_1}{\tau_1},\ldots,\frac{W_K^\ast-\mu_K}{\tau_K}
\right)
\]
where the $\mu_i$ and $\tau_i$ are the conditional means and standard deviations
of the permutation distribution of $W_i^\ast$, given the tie pattern
${\cal T}$. The test is performed conditionally given 
${\cal T}$. See \cite{Lehmann06}, Chapter 2,
for a discussion of the issues involved
in the case of the Wilcoxon test. 
For equal sample sizes across treatments the test based on $\hat{S}_{\max}$ is equivalent to the one based on $S$.

When treatment effects are expected to lower scores compared to those under control one would
of course reject $H$ for low values of
\[
\hat{S}_{\min} =\min\left(
\frac{W_1^\ast-\mu_1}{\tau_1},\ldots,\frac{W_K^\ast-\mu_K}{\tau_K}
\right)
\]
For the two-sided test we should reject $H$ when
\[
\hat{S}_{\rm abs} =\max\left(
\frac{|W_1^\ast-\mu_1|}{\tau_1},\ldots,\frac{|W_K^\ast-\mu_K|}{\tau_K}
\right)
\]
is too large, or equivalently when either $\hat{S}_{\min}$ is too small or $\hat{S}_{\max}$ is too large, with rejection limits symmetric around zero.

\section{Conditional Moments}
Here we focus on the issue of ties and work out the details for using the multivariate normal
distribution as a large sample approximation for the joint distribution of $(W_1^\ast,\ldots,W_K^\ast)$ under
the randomization model, assuming that $H$ holds.

Following \cite{Miller81} we will use the equivalent Mann-Whitney form of the
rank sum statistic, i.e.,
$
W_i^\ast = W_{X_0,X_i}^\ast+n_i(n_i+1)/2
$
where for samples $X_0=(X_{01},\ldots,X_{0n_0})$ and
$X_i=(X_{i1},\ldots,X_{in_i})$ the Mann-Whitney statistic in case of
ties is defined by
\[
W_{X_0,X_i}^\ast=\sum_{j=1}^{n_0}\sum_{k=1}^{n_i} I_{\{X_{0j}<X_{ik}\}}+
\frac{1}{2}\sum_{j=1}^{n_0}\sum_{k=1}^{n_i} I_{\{X_{0j}=X_{ik}\}}
\]
The Mann-Whitney form avoids having to re-rank observations for
each separate pair of samples. 

It will suffice to get the means, variances and covariances of
$W_{X_0,X_1}^\ast$ and
$W_{X_0,X_2}^\ast$.
For ease of exposition, we refer from now on to $X_0=(X_{01},\ldots,X_{0n_0})$
as $X=(X_1,\ldots,X_{n_0})$, to $X_1=(X_{11},\ldots,X_{1n_1})$ as
$Y=(Y_1,\ldots,Y_{n_1})$ and to
$X_2=(X_{21},\ldots,X_{2n_2})$ as
$Z=(Z_1,\ldots,Z_{n_2})$, where the three samples $X$, $Y$, and $Z$
are randomly selected without replacement from the $N$ (midrank)
observations
$w_1,\ldots, w_N$ with the tie pattern ${\cal T}$. This randomization model 
forms the basis for all probabilities, expectations, variances and covariances.

First we obtain the expectation of $W_{X,Y}^\ast$
under the randomization distribution.
We condition on the 
sampled tie pattern ${\cal T}$ and not on the tie pattern 
of just the $n_0+n_1$ observations. Let 
$D_i = \sum_{\nu=1}^i d_\nu, i=1,\ldots, e$, then
\[
E(W_{X,Y}^\ast)=\sum_{j=1}^{n_0}\sum_{k=1}^{n_1} P(X_j<Y_k)+\frac{1}{2}\sum_{j=1}^{n_0}\sum_{k=1}^{n_1} P(X_j=Y_k)
\]
where
\[
P(X_j<Y_k)=\sum_{\nu=1}^eP(X_j<Y_k|Y_k=u_\nu)P(Y_k=u_\nu)
=\sum_{\nu=1}^e\frac{d_\nu(D_\nu-d_\nu)}{N(N-1)}
\]
and
\[
P(X_j=Y_k)=\sum_{\nu=1}^eP(X_j=Y_k|Y_k=u_\nu)P(Y_k=u_\nu)
=\sum_{\nu=1}^e\frac{d_\nu(d_\nu-1)}{N(N-1)}
\] 
Thus
\[
E(W_{X,Y}^\ast)=\frac{n_0n_1}{N(N-1)}\sum_{\nu=1}^e \left(d_\nu(D_\nu-d_\nu)+\frac{1}{2}d_\nu(d_\nu-1)\right)
=\frac{n_0n_1}{2}
\]

For the variance of $W_{X,Y}^\ast$ we need a generalization of formula (1.35) in \cite{Lehmann06}, namely
\begin{equation}
\var(W_{X,Y}^\ast)=\frac{n_0n_1(n_0+n_1+1)}{12}-\frac{n_0n_1\sum_{\nu=1}^e(d_\nu^3-d_\nu)}{12(n_0+n_1)(n_0+n_1-1)}
\label{varties}
\end{equation}
which assumes $n_0+n_1=N$ in the underlying randomization model.

While it is possible to follow the approach used in obtaining $E(W_{X,Y}^\ast)$ also for $\var(W_{X,Y}^\ast)$,
it is extremely laborious and not conducive to obtaining succinct expressions.
The approach taken here is based on the following well known identity
\begin{equation}
\var({\cal Q})=E\left(\var(\cal{Q}|\cal{U})\right)+
\var\left(E(\cal{Q}|\cal{U})\right)
\label{var.cond}
\end{equation}
where $\cal{Q}$ is any random variable for which the second moment exists
and $\cal{U}$ can represent any random vector. In our situation we will substitute the 
Mann-Whitney statistic $W_{X,Y}^\ast$ for $\cal{Q}$ in (\ref{var.cond}).
Computing $\var(W_{X,Y}^\ast)$ under the permutation distribution, we take 
for $\cal{U}$ the random sample of size $N-n_0-n_1=\tilde{n}$ that is left over when splitting the 
total set of $N=d_1+\ldots+d_e$ midranks into the three samples $X$, $Y$ and $\cal{U}$
of respective sizes $n_0, n_1$, and $\tilde{n}$.  The random sample 
$\cal{U}$ can be characterized by the 
counts $G_1,\ldots, G_e$ with which the unique midrank
values $u_1,\ldots,u_e$ appear 
in the sample $\cal{U}$. We now view the samples $X$ and $Y$ as any random split
of the  $u_1,\ldots, u_e$, now with remaining multiplicities
$d'_\nu=d_\nu-g_\nu$, $\nu=1,\ldots,e$, where $g_\nu$ represents the realized count of $G_\nu$.
Note that we have
\[
N'=\sum_{\nu=1}^e d'_\nu= \sum_{\nu=1}^e d_\nu+\sum_{\nu=1}^e g_\nu=N-\tilde{n}=n_0+n_1
\]
Applying the variance formula
(\ref{varties}) in this conditional permutation situation we get
\[
\var\left(
W_{X,Y}^\ast|\cal{U}
\right)=\frac{n_0n_1(n_0+n_1+1)}{12}-\frac{n_0n_1}{12(n_0+n_1)(n_0+n_1-1)}\sum_{\nu=1}^e({d'_\nu}^3-d_\nu')
\]
Since  $E(W_{X,Y}^\ast|{\cal{U}})=n_0 n_1/2$ is constant in $\cal{U}$ we have $\var\left(E(W_{X,Y}^\ast|\cal{U})\right)=0$. We need
to find the expectation of 
\[
\sum_{\nu=1}^e({d'_\nu}^3-d_\nu')=\sum_{\nu=1}^e((d_\nu-G_\nu)^3-(d_\nu-G_\nu))
\]
and use it in the above formula to get the permutation variance of $W_{X,Y}^\ast$. It suffices to find
the expectations
\[
E\left((d_\nu-G_\nu)^3\right)\qquad \mbox{and}\qquad E\left(d_\nu-G_\nu\right)
\] 
where $G_\nu$ is simply a hypergeometric random variable, counting the number of 
times the unique value $u_\nu$ is included in the sample $\cal{U}$ when sampling $\tilde{n}$ times 
without replacement from the finite population of $N$ midrank values with the given tie pattern ${\cal T}$. This
expectation calculation
is straighforward when using
the formulas
\[
E\left(G_\nu(G_\nu-1)(G_\nu-2)\right)=\frac{\tilde{n}(\tilde{n}-1)(\tilde{n}-2)d_\nu(d_\nu-1)(d_\nu-2)}{N(N-1)(N-2)}
\]
and
\[
E\left(G_\nu(G_\nu-1)\right)=\frac{\tilde{n}(\tilde{n}-1)d_\nu(d_\nu-1)}{N(N-1)}
\]
After some algebra with a view toward factorization one finds
\[
E\left((d_\nu-G_\nu)^3\right)-E\left(d_\nu-G_\nu\right)\hspace{4in}
\]
\begin{eqnarray*}
&=&
\qquad \frac{(N-\tilde{n})(N-\tilde{n}-1)(N-\tilde{n}-2)}{N(N-1)(N-2)} d_\nu(d_\nu-1)(d_\nu+1)\\
&&\hspace{2in}+ \frac{3(N-\tilde{n})(N-\tilde{n}-1)\tilde{n}}{N(N-1)(N-2)} d_\nu(d_\nu-1)
\end{eqnarray*}
and  thus
\begin{eqnarray}
\var(W_{X,Y}^\ast)&=& 
\frac{n_0 n_1(n_0+n_1+1)}{12}\nonumber \\ 
&&\qquad -
\frac{n_0n_1(n_0+n_1-2)}{12N(N-1)(N-2)}\sum_{\nu=1}^e d_\nu(d_\nu-1)(d_\nu+1)\nonumber \\
&& \qquad -\s \frac{n_0n_1(N-n_0-n_1)}{4N(N-1)(N-2)}\sum_{\nu=1}^e d_\nu(d_\nu-1)\label{varformula}
\end{eqnarray}
Note that (\ref{varformula}) reduces to (\ref{varties}) when $N=n_0+n_1$.

To compute $\cov(W_{X,Y}^\ast,W_{X,Z}^\ast)$ we note that 
$W_{X,Y}^\ast+W_{X,Z}^\ast=W_{X,(Y,Z)}^\ast$, where
$(Y,Z)$ is the concatenated vector of the $Y$ and $Z$ samples.
From that we obtain
\[
\cov(W_{X,Y}^\ast,W_{X,Z}^\ast)=\frac{1}{2}\left(\var\left(W_{X,(Y,Z)}^\ast\right)-\var(W_{X,Y}^\ast)-\var(W_{X,Z}^\ast)\right)
\]
where the variances on the right side can all be expressed using (\ref{varformula}) with appropriate sample sizes. This leads
to
\begin{equation}
\cov\left(W_{X,Y}^\ast,W_{X,Z}^\ast\right)=\frac{n_0n_1n_2}{12}-\frac{n_0n_1n_2}{12N(N-1)(N-2)}
\sum_{\nu=1}^e d_\nu(d_\nu-1)(d_\nu-2)
\label{covformula}
\end{equation}
The tie correction terms in (\ref{varformula}) and 
(\ref{covformula}) disappear when there are no ties, i.e., when
$d_i=1$ for all $i$. In the presence of ties the relative magnitude
of these terms is small, except for extreme tie patterns, as illustrated
in Section~\ref{Approximate}.

\section{Approximate Significance Probabilities\label{Approximate}}
The covariance structure of $W_{X_0,X_i}^\ast, i=1,\ldots,K$, allows us to view these random variables
as having the same covariance structure as the simple linear expressions $n_i V_0+V_i, i=1,\ldots,K$,
in terms 
of the independent random variables
$V_0, V_1, \ldots, V_K$. With the variances of the latter denoted by $\var(V_i)=\sigma_i^2, i=0,\ldots,K$, we have  for $i,j\neq 0$
$\cov(n_i V_0+V_i, n_j V_0+V_j)=n_in_j \sigma_0^2$
for $i \neq j$ and
$\tau_i^2=\var(n_i V_0+V_i)=n_i^2\sigma_0^2+\sigma_i^2$,
i.e., $\sigma_i^2=\tau_i^2-n_i^2 \sigma_0^2$. Matching appropriate terms in (\ref{varformula}) and
(\ref{covformula}) with $\tau_i^2$ and $\sigma_0^2$ we get
\[
\sigma_i^2=\frac{n_0n_i}{12}\left(
n_0+1-\frac{3}{N(N-1)}\sum_{\nu=1}^e d_\nu(d_\nu-1)\right.\hspace{2in}
\]
\[
\left.\hspace{2in}-\frac{n_0-2}{N(N-1)(N-2)}\sum_{\nu=1}^e d_\nu(d_\nu-1)(d_\nu-2)
\right)
\]
for $i=1\ldots, k$ and\[
\sigma_0^2=\frac{n_0}{12}-\frac{n_0}{12N(N-1)(N-2)}
\sum_{\nu=1}^e d_\nu(d_\nu-1)(d_\nu-2)
\]
When appealing to the asymptotic multivariate normality
of the statistics $W_{X_0,X_1}^\ast, \ldots, W_{X_0,X_K}^\ast$
for large $n_0, n_1,\ldots, n_K$ (to be shown in the Appendix),
we may add the normality assumption to the independent terms $V_0, V_1, \ldots, V_K$ in
the above asymptotic representation of $W_{X_0,X_i}^\ast$ in terms of $n_i V_0+V_i$.
Since we will use the $W_{X_0,X_i}^\ast$ in standardized form in the Steel statistic we may also assume 
that $V_0, V_1,\ldots,V_K$ all have mean zero.
This
reduces the approximate upper tail conditional probability calculation to
\[
P\left(\max_{i=1,\ldots, K}\left(
\frac{W_{X_0,X_i}^\ast-n_0n_i/2}{\tau_i}\right) \geq U \right) \approx 
1-P\left(\max_{i=1,\ldots, K}\left(
\frac{n_i V_0+V_i}{\tau_i}\right) \leq U \right)
\]
\begin{equation}
= 1-\int_{-\infty}^\infty \prod_{i=1}^K \Phi\left(
\frac{U\tau_i-n_i\sigma_0 z}{\sigma_i}\right) \varphi(z)dz
\label{upper}
\end{equation}
where $\Phi(z)$ and $\varphi(z)$ are the standard normal
distribution function and density, respectively.
For the one-sided test in the other direction the calculation is
\[
P\left(\min_{i=1,\ldots, K}\left(
\frac{W_{X_0,X_i}^\ast-n_0n_i/2}{\tau_i}\right) \leq L \right) \approx 
1-P\left(\min_{i=1,\ldots, K}\left(
\frac{n_i V_0+V_i}{\tau_i}\right) \geq L \right)
\]
\begin{equation}
= 1-\int_{-\infty}^\infty \prod_{i=1}^K \left[1-\Phi\left(
\frac{L\tau_i-n_i\sigma_0z}{\sigma_i}\right)\right] \varphi(z)dz
\label{lower}
\end{equation}
while for the two-sided test it becomes
\[
P\left(\max_{i=1,\ldots, K}\left(
\frac{|W_{X_0,X_i}^\ast-n_0n_i/2|}{\tau_i}\right) \geq U \right) \approx 
1-P\left(\max_{i=1,\ldots, K}\left(
\frac{|n_i V_0+V_i|}{\tau_i}\right) \leq U \right)
\]
\begin{eqnarray*}
= 1-\int_{-\infty}^\infty \prod_{i=1}^K \left[\Phi\left(
\frac{U\tau_i-n_i\sigma_0 z}{\sigma_i}\right)- \Phi\left(
\frac{-U\tau_i-n_i\sigma_0 z}{\sigma_i}\right)\right]\varphi(z)dz
\end{eqnarray*}

The Appendix shows that asymptotic multivariate normality
of $W_{X_0,X_1}^\ast, \ldots, W_{X_0,X_K}^\ast$
holds for the conditional randomization context,
using conditional means, variances
and covariances given the tie pattern ${\cal T}$. The conditions for 
this are:\\
$C_1$:
\quad $n_0,n_1,\ldots, n_K\rightarrow \infty$ with 
$n_i/N\rightarrow \rho_i>0$, $i=0,\ldots,K$, and \\
$C_2$:
\quad $\max(d_i/N)\leq 1-\epsilon$ as $N\rightarrow \infty$ for some 
$\epsilon >0$.

The effectiveness of this approximation is easily assessed through simulation. 
Figure~\ref{fig1} shows an example
for the simulated distribution of $\hat{S}_{\rm max}$ 
in the case of 3 samples 
with $n_0=n_1=n_2=100$ from continuous distributions (standard normal here, although it is
not relevant here),
comparing its upper tail probabilities with its asymptotic
approximation (\ref{upper}). 
Figure~\ref{fig2} shows the same comparison after rounding
the 300 sample values to a single decimal place, resulting 
in 51 unique values. Here the sampled continuous distribution
is relevant, since the rounding will
result in different tie patterns for different sampled distributions.
In both displayed cases the simulated distribution is based on
${\sf Nsim=1e+06}$ random splits of the 300 (rounded) observations
into samples of size 100 each, thus approximating the conditional permutation distribution.
The approximations are 
quite satisfactory. Some more examples are presented in 
Section~\ref{AppQual}. The reader is invited to further explore this issue
with the {\sf \cite{R}} package {\tt kSamples}, see \cite{Scholz16}.

\begin{figure}[htb]
\begin{center}
\includegraphics[width=3.3in]{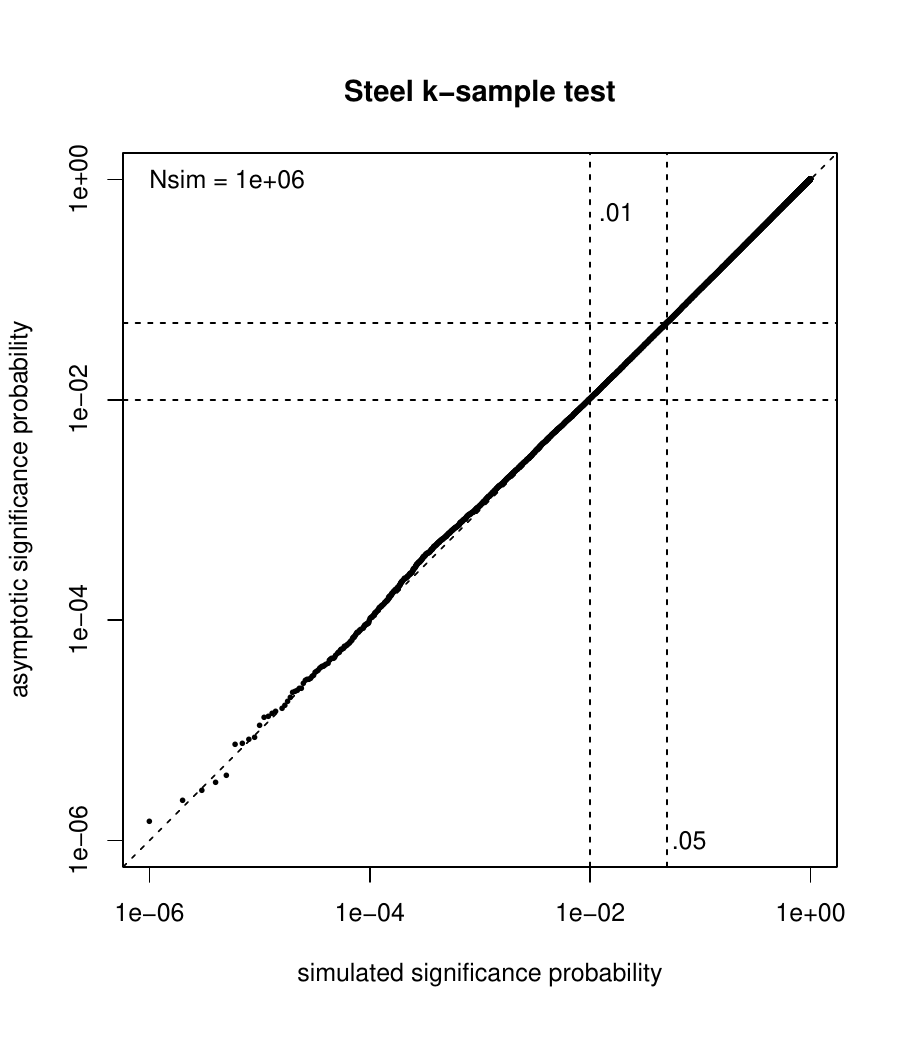}
\caption{Approximation Quality for (\ref{upper}),
$K=2$, $n_0=n_1=n_2=100$
\label{fig1}}
\end{center}
\end{figure}

\begin{figure}[htb]
\begin{center}
\includegraphics[width=3.3in]{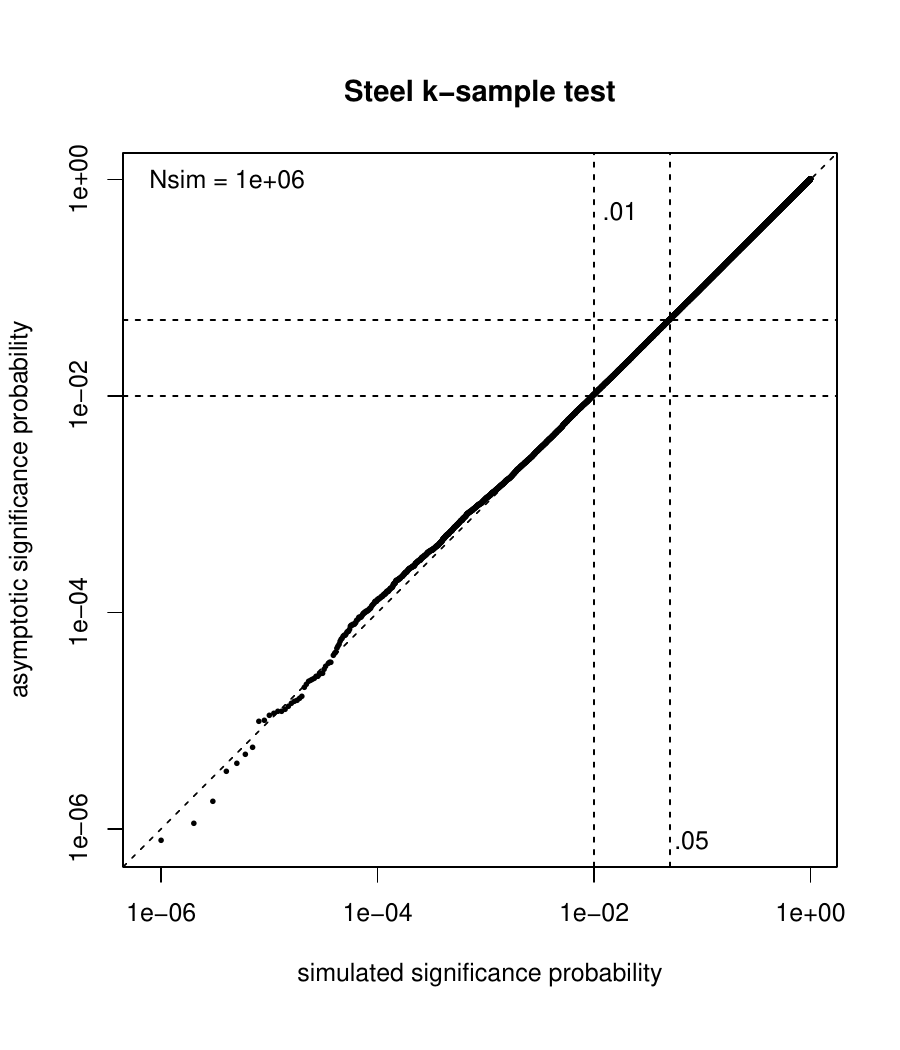}
\caption{Approximation Quality for (\ref{upper}) after Rounding,
$K=2$, $n_0=n_1=n_2=100$
\label{fig2}}

\end{center}
\end{figure}

\section{An Example}
As a computational validation 
example we use the data of the effect of birth condition on the IQ of 24 
girls, as presented in \cite{Steel59}. 
Since the ``treatment'' effect would presumably
lead to possibly lower IQ values we take the minimum version $\hat{S}_{\min}$ as test
statistic. Its value is
$-1.7713$ with simulated significance probability of
$0.10474$ based on 100000 simulated splits of the 24 observed IQ's
into equal groups of 6 each. The normal approximation based on (\ref{lower})
gives $0.0946$, reasonably accurate considering the smallness
of the involved sample sizes, 6 and 6 for each comparison.
The observed values of $W_{X_0,X_i}^\ast$ were $7, 17, 12.5$, respectively.
The values used for standardizing the $W_{X_0,X_i}^\ast$ and implementing
the normal approximation are $\sigma_0=0.7062328$, $\mu_i=18$,
$\sigma_i=4.540007$, $i=1,2,3$, and $\tau=6.210249$.

\section{Approximation Quality\label{AppQual}}

We offer here two types of comparison.
In the first we compare the approximation quality of the Steel test based
on $\hat{S}_{\max}$ with the 
normal approximation of the 2-sample Wilcoxon test based on the control sample
and the first treatment sample. This is illustrated in Figure~\ref{fig3} and shows that 
the normal approximation in the 2-sample case seems slightly better. We believe that this
is due to the smaller proportion of ties in the first two samples alone.  

As illustrated by \cite{vdWiel02}, exact computations can become prohibitively time consuming 
once sample 
sizes or the number of treatments become large. Power considerations may 
lead to larger sample sizes. In that case
even simulations become time consuming and asymptotic approximations as presented here
offer a reasonable practical alternative.
\begin{figure}
\begin{center}
\includegraphics[width=4in]{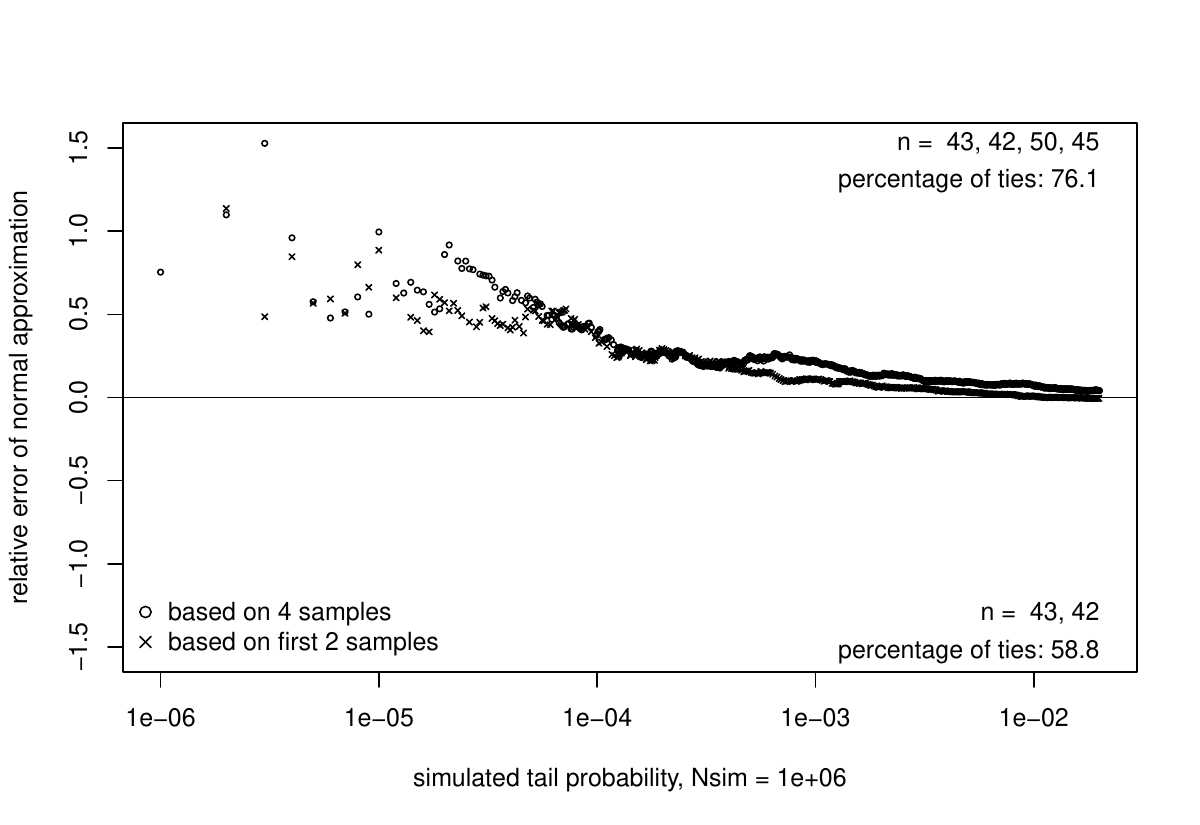}
\caption{Comparing Normal Approximation Quality of the 
Steel Test
with Two-sample Wilcoxon Test\label{fig3}}
\end{center}
\end{figure}

In the second comparison we examine what would happen if we ignore the presence of ties
and thus use the wrong standard deviations $\tau_1,\ldots, \tau_K$ in $\hat{S}_{\max}$
and also the wrong $\tau_i, \sigma_i$, and $\sigma_0$ in the normal approximation and compare that
with the proper treatment. This is done in the extreme case when there are only two distinct sample
values, with roughly equal frequency across the samples, and then again when there are only 9
distinct values. The results are displayed in Figs.~\ref{fig4} and \ref{fig5}, respectively.
In Fig.~\ref{fig5} there seems to be little difference between the two actions. Examining the 
expressions (\ref{varformula}) and (\ref{covformula}) for the variances and covariances it becomes
clear that the ties contribute only small adjustments to the same expressions when there are no ties 
($d_\nu=1$ for all $\nu$), except when the tie pattern is extreme, as illustrated 
in Fig.~\ref{fig4}.

\begin{figure}
\begin{center}
\includegraphics[width=3.5in]{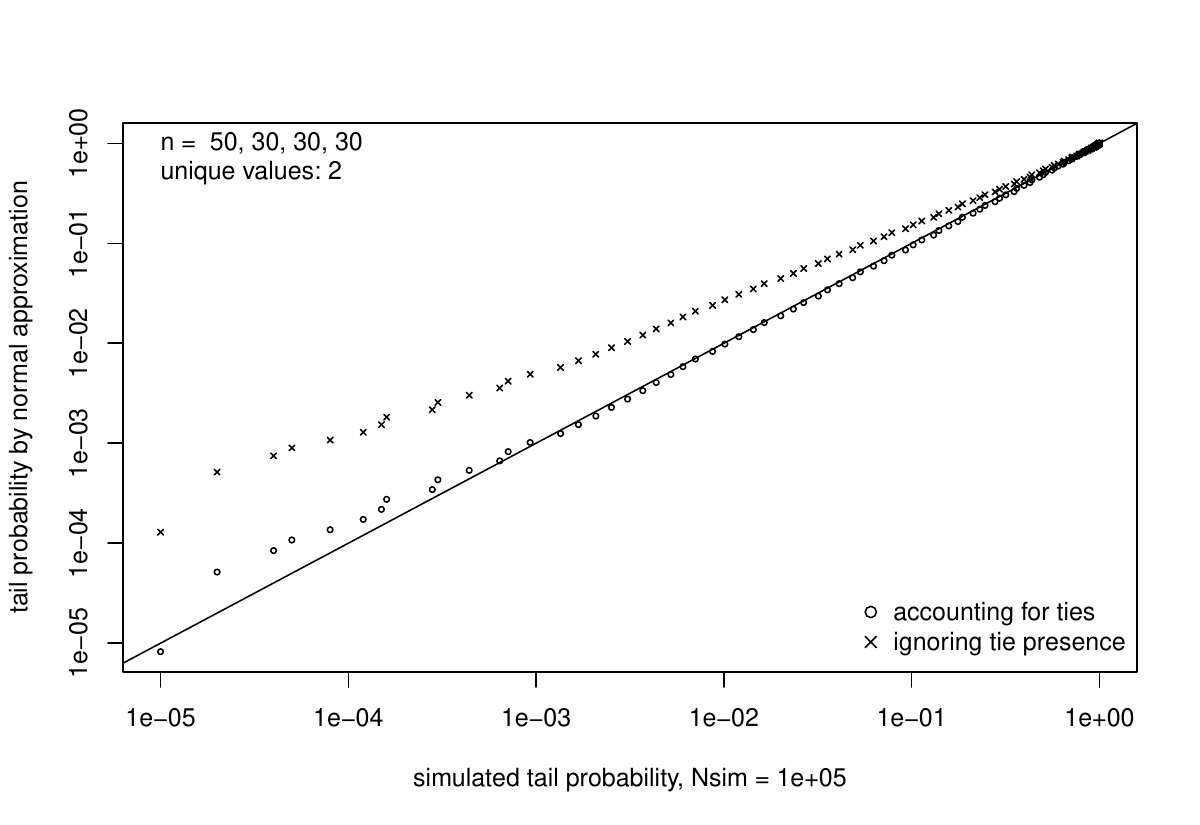}
\caption{Normal Approximation Quality of Steel Test,
with and without Adjusting Extreme Ties \label{fig4}}
\end{center}
\end{figure}

\begin{figure}
\begin{center}
\includegraphics[width=3.5in]{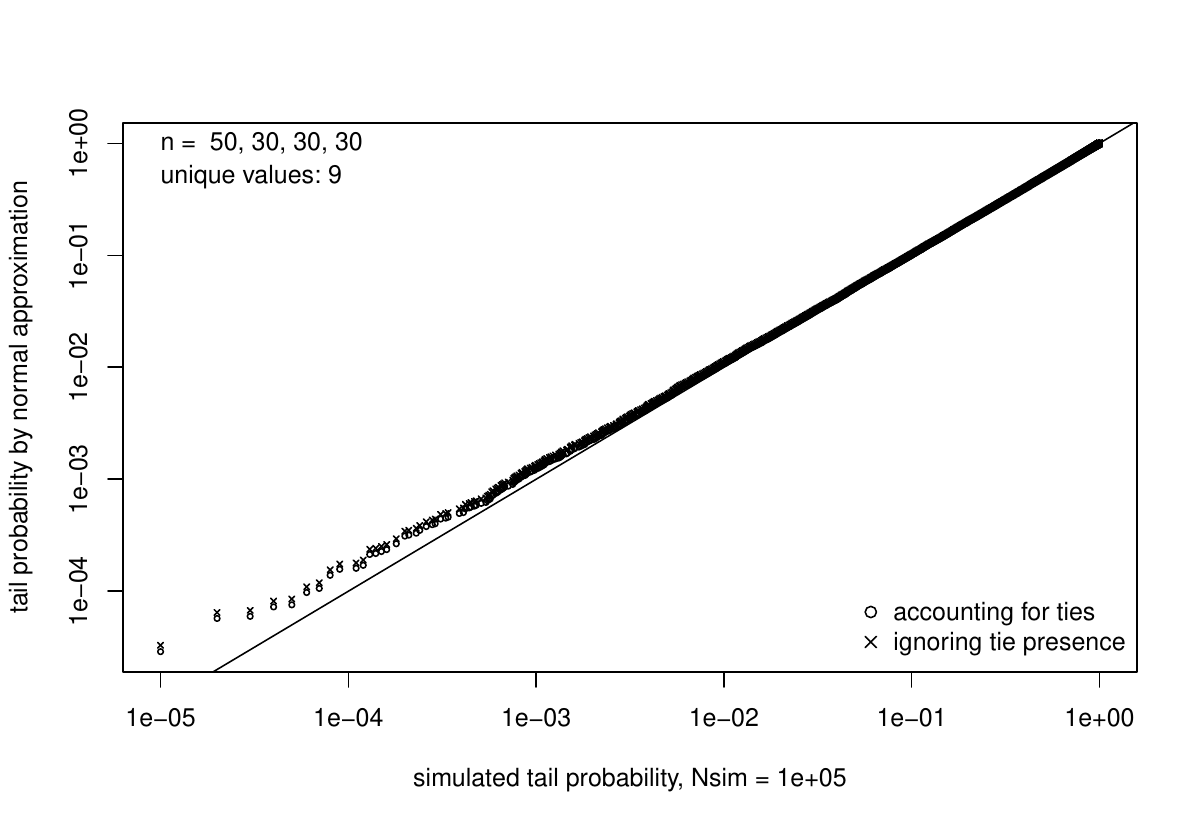}
\caption{Normal Approximation Quality of Steel Test,
with and without  Adjusting for Many Ties\label{fig5}}
\end{center}
\end{figure}

\section{Simultaneous Confidence Bounds}

Here we assume that the control and treatment samples
are independent random samples from populations
with continuous cumulative distribution functions
$F_0(x)$ and $F_i(x)=F_0(x-\Delta_i)$, $i=1,\ldots,K$,
respectively.
The main purpose of the continuity assumption is 
to make the above shift model meaningful for any
values of $\Delta_i$. Our goal is to turn Steel's simultaneous test procedure into
simultaneous confidence bounds or intervals
for $\Delta_1,\ldots,\Delta_K$, as suggested at the 
end of section 5A in \cite{Lehmann06}.

Let $D_{i(j)}$ denote the $j$-th ordered value (ascending order) of the $n_0n_i$
paired differences between the observations in treatment sample $i$ and those of the control
sample. 
By simple extension of results in \cite{Lehmann06}, 
pages 87 and 92, the following equations hold, 
relating the null distribution of the 
Mann-Whitney statistics
$W_{X_0,X_i}$ and the coverage probabilities of the $D_{i(j_i)}$ for any set of
$j_1,\ldots,j_K$ with $1\le j_i \le n_0 n_i$.
\begin{equation}
P_\Delta(\Delta_i \le D_{i(j_i)}, i=1,\ldots,K)=
P_0(W_{X_0,X_i}\le j_i -1, i=1,\ldots,K)
\label{bound1}
\end{equation}
and
\begin{equation}
P_\Delta(\Delta_i \ge D_{i(j_i)}, i=1,\ldots,K)=
P_0(W_{X_0,X_i}\le n_1 n_i -j_i, i=1,\ldots,K)
\label{bound2}
\end{equation}
where $P_\Delta$ refers to the distribution under $\Delta=(\Delta_1,\ldots,\Delta_K)$
and $P_0$ refers to the joint null distribution of the $W_{X_0,X_i}$ when all sampled 
distributions are the same and continuous. There are $K$ indices $j_i$ that can be manipulated
to affect the achieved confidence level. We discuss only the 
treatment of equation (\ref{bound1}). The case of
(\ref{bound2}) is treated analogously.
To limit the computational complexity
standardized versions of the $W_{X_0,X_i}$, i.e.,  $(W_{X_0,X_i}-n_0n_i/2)/\tau_i$ 
are used to choose a common value for $(j_i -1-n_0n_i/2)/\tau_i$ from the multivariate normal approximation 
for the $W_{X_0,X_i}$ (see \cite{Miller81}), and reduce that 
to integer values for $j_i$, rounding up, rounding down, and rounding to the nearest integer. These sets of three
integers  $j_i$ are then used in approximating the actual probabilities
$P_0(W_{X_0,X_i}\le j_i -1, i=1,\ldots,K)$, and from these three coverage probabilities 
the one that is closest to the nominal confidence level $\gamma$ and $\ge \gamma$
and also the one that is closest without the restriction $\ge \gamma$ are chosen.
This is implemented in the {\sf R} package {\tt kSamples}.

Confidence intervals are constructed by using upper and lower confidence bounds, each with
same confidence level of $(1+\gamma)/2$.

When the original sample data appear to be rounded, and especially when there are ties,
one should widen the computed intervals or bounds by the rounding $\epsilon$, as illustrated
in \cite{Lehmann06}, pages 85 and 94. For example, when all sample values appear to end in one of $.0, .2, .4, .6, .8$,
the rounding $\epsilon$ would be $.2$. Ultimately, this is a judgment call for the user. Such widening
of intervals will make the actually achieved confidence level $\ge$ the stated achieved level. 

\section{All Pairwise Comparisons}
While Steel's procedure, as discussed so far, focuses on pairwise comparisons
between a single control sample and several treatment samples, the conditional moment
results derived here also apply to the conditional distributions of all pairwise
Wilcoxon comparisons involving $K$ treatment samples. Such all pairwise tests were
proposed by \cite{Dwass60} and
\cite{Steel60} and were more recently examined by \cite{Neu01}.
The resulting 
approximate computation of significance probabilities, when dealing with the maximum 
or minimum of all standardized test statistics, can be achieved using a multivariate
normal distribution with dimension ${K \choose 2}$, using the derived conditional moment structure.
To this we add the following:
$\cov\left(W_{X_1,X_2}^\ast,W_{X_3,X_1}^\ast\right)=-\cov\left(W_{X_1,X_2}^\ast,W_{X_1,X_3}^\ast\right)$
and
$\cov\left(W_{X_1,X_2}^\ast,W_{X_3,X_4}^\ast\right)=0$,
all understood to be conditional on the tie pattern. The first follows from 
$W_{X_3,X_1}^\ast+W_{X_1,X_3}^\ast=n_1n_2$ and the second from a covariance analog to
(\ref{var.cond}), conditioning on the pooled samples of sizes $n_1+n_2$ and $n_3+n_4$
and realizing that their further splits into $(n_1,n_2)$ and $(n_3,n_4)$ by randomization are independent of each other. Asymptotic multivariate normality should hold under the same conditions 
$C_1$ and $C_2$ (omitting index zero) and same line of approach 
as presented in the Appendix.

\section*{Appendix}
Under sampling from distributions
$F_0=F_1=\ldots=F_K$ the asymptotic multivariate normality is shown 
in \cite{Lehmann63} (as pointed out by \cite{Miller81}),
and again by \cite{MunzelHothorn01} who address
the widened scope of sampling from discontinuous
distributions as well, with possibly unequal treatment sample sizes.

Here we complement these results by showing asymptotic multivariate normality of $(W_{X_0,X_1}^\ast,\ldots,  W_{X_0,X_K}^\ast)$ under the treatment randomization model, assuming conditions $C_1$ and $C_2$ hold,
see end of Section~\ref{Approximate}.

We have $N$ subjects split randomly into groups of sizes
$n_0,n_1,\ldots, n_K$ with $n_0+n_1+\ldots+n_K=N$.
To these subjects
we respectively administer either a control (placebo) 
or one of treatments $1$ through $K$. Under the hypothesis of no treatment effect (i.e., 
treatments and control do not affect the subjects in differential ways) we can assume
that any scores that are measured as experimental outcome for each subject are preordained.
We allow for the possibility of ties. Since we want to employ ranks for testing purposes
we will work with the vector of midranks of all $N$ subjects, denoted by 
$v_N=(v_{N1},\ldots,v_{NN})$ with average $\bar{v}_{N}=\sum_i v_{Ni}/N
=(N+1)/2$. Without loss 
of generality we assume that labeling of the subjects is such that 
$v_{N1}\leq v_{N2}\leq \ldots \leq v_{NN}$. By $u_{N1}< \ldots < u_{Ne}$ we denote their
unique values with respective multiplicities $d_1,\ldots, d_e\geq 1$. We should also write
$d_{Ni}$ and $e_N$ instead of $d_i$ and $e$, but for notational simplicity we refrain
from doing so.

We employ  i.i.d.~uniform $U(0,1)$ 
random variables
$U_1,\ldots,U_N$ with order statistics $U_{(1)} < \ldots < U_{(N)}$ and denote
by $R_i$ the rank of $U_i$ among $U_1,\ldots,U_N$. With probability one
there are no ties among these
$R_i$. 
We define the following indicator random variables to prescribe 
which subject will belong to the control or which treatment group.
For $\ell=0,1,\ldots,K$ and using $m_{\ell-1}=n_0+n_1+\ldots+n_{\ell-1}$ with $m_0=n_0$ and $m_{-1}=0$
we define
\[
J_{\ell i}=\left\{\begin{array}{lll}
1 & \mbox{if} & R_i\in \{m_{\ell-1}+1,\ldots, m_{\ell-1}+n_\ell\}=A_\ell\\
0 & \mbox{else} &
\end{array}
\right.
\]
with the understanding that midrank $v_{Ni}$ is associated with the 
control or $\ell^{\rm th}$ treatment group whenever $J_{0,i}=1$ or
$J_{\ell,i}=1$, respectively.  We denote the resulting samples
of allocated $v_{Ni}$ values
by $X, Y_1,\ldots, Y_K$, respectively.

Under this randomization scheme all 
${N\choose n_0,n_1,\ldots,n_K}$ possible partitions of the $N$ subjects
into groups of respective sizes $n_0, n_1,
\ldots, n_K$ are equally likely. This is the basis for the null distribution 
of any statistic based symmetrically on the scores within each group. 

Comparing each treatment against the same control we form separate 
Mann-Whitney statistics, i.e., for $\ell=1,\ldots,K$ using $X$ as control
and $Y_\ell$ as treatment sample we define
\begin{eqnarray*}
W_{X,Y_\ell}^\star&=&
\sum_{1\leq i \neq j \leq N} J_{0i}J_{\ell j}\left[I_{[v_{Ni}<v_{Nj}]}
+\frac{1}{2}I_{[v_{Ni}=v_{Nj}]}
\right]
= \sum_{1\leq i \neq j \leq N} J_{0i}J_{\ell j}\psi_{ij}(v_N)\\
&&\quad
\mbox{with}\quad \psi_{ij}(v_N)=\left[I_{[v_{Ni}<v_{Nj}]}
+\frac{1}{2} I_{[v_{Ni}=v_{Nj}]}
\right]=1-\psi_{ji}(v_N)
\end{eqnarray*}
Here the $i\neq j$ under the summation results from 
$J_{0i}J_{\ell i}=0$ for $\ell\neq 0$.

From now on we use the convention that
$\sum_{i,j}, \sum_{i,j,r},  \sum_{i,j,r,s}$ indicate summations with indices running
from $1$ to $N$ and all shown indices are distinct, unless explicitly indicated 
otherwise.

$\psi_{ij}+\psi_{ji}=1$ implies
$\bar{\psi}(v_N)=\frac{1}{N(N-1)}\sum_{i,j} 
\psi_{ij}(v_N)=1/2$ and $\psi_{ii}=1/2$.
\[
W_{X,Y_\ell}^\star=\sum_{i,j} J_{0i}J_{\ell j}(\psi_{ij}(v_N)-\bar{\psi}(v_N))+
n_0n_\ell \bar{\psi}(v_N)
=\sum_{i,j} J_{0i}J_{\ell j}\hat\psi_{ij}(v_N)+
n_0n_\ell/2
\]
where
\begin{equation}
\hat\psi_{ij}(v_N)=\psi_{ij}(v_N)-\bar{\psi}(v_N)
\quad \mbox{with}\quad
\sum_{i,j} \hat\psi_{ij}(v_N)=0\quad
\mbox{and} \quad \hat{\psi}_{ii}=0
\label{psi0}
\end{equation}
Note that the ordered labeling of the midranks $v_{Ni}$ yields
\begin{eqnarray}
\hat{\psi}_{\cdot n}(v_{N})&=&
\sum_{i=1}^N\hat{\psi}_{in}(v_{N})=\sum_{i=1}^N\psi_{in}-\frac{N}{2}
\nonumber\\
&=&d_1+\ldots+d_{n-1}+\frac{d_n+1}{2}-\frac{N+1}{2}
=v_{Nn}-\bar{v}_N
\label{psi1}
\end{eqnarray}
and
\begin{equation}
\hat{\psi}_{n\cdot}(v_{N})=
\sum_{j=1}^N\hat{\psi}_{nj}(v_{N})=-(v_{Nn}-\bar{v}_N)
\label{psi2}
\end{equation}

The corresponding treatment rank sums
\[
W_{n_0n_\ell}^\star=W_{X,Y_\ell}^\star+\frac{n_\ell(n_\ell+1)}{2}
\]
amount to re-ranking the midranks within each 
pair of groups (control group and $\ell^{\rm th}$ treatment group)
and then take the sum of (modified) treatment midranks. While the individual statistics
may well be approximately normal (still needs proof in our randomization context),
they will not be independent, since each
comparison statistic involves
the same control group. 
Here we aim to show that a $K$-variate multivariate
normal distribution is an appropriate approximation under conditions
$C_1$ and $C_2$. Our 
proof
follows closely the approach of \cite{Lehmann06} which was based on \cite{Hajek61}.

Without the re-ranking the joint asymptotic normal distribution of the corresponding 
treatment rank sums in the randomization setting
is already established at the end of \cite{Hajek61}. 
The relative merits of these two approaches are not clear
and it may be worth investigation.

Let
\[
T_\ell=W_{X,Y_\ell}^\ast-\frac{n_0n_\ell}{2}=
\sum_{i,j}J_{0i}J_{\ell j}\hat{\psi}_{ij}(v_N)\quad
\mbox{with \quad $E(T_\ell)=0$}
\]

H\'ajek's projection method, see \cite{Lehmann06}  p.~362ff,
approximates $T$ by 
a sum of iid random variables for which asymptotic normality is
easily established.
This approximation is the projection of $T$ onto the space of functions of
form $\sum_{r=1}^N f(U_r)$ and is obtained by taking 
\[
f(U_r)=E(T_\ell|U_r)=\sum_{i,j} E(J_{0i}J_{\ell j}|U_r)\hat{\psi}_{ij}(v_N)
\]
To evaluate this we need to find
the conditional probabilities
\[
g_{mnijr}(u)=P(R_i=m,R_j=n|U_r=u)
\]
for $i\neq j$ and $m\neq n$.

Because of the exchangeability of the $U_r$'s we get
that for distinct $i,j,r$ the
$g_{mnijr}(u)=g_{mn}(u)$ are
independent of $i,j,r$,
and thus
\[
\sum_{i\neq j=1;i,j\neq r}^N P(R_i=m,R_j=n|U_r=u)=
(N-1)(N-2)g_{mn}(u)
\]
By exchangeability we also get
$
g_{mnrjr}(u)=h_{mn}(u)$
and 
$g_{mnirr}(u)=h_{nm}(u)$ independent of $i,j\neq r$.
\begin{eqnarray*}
\sum_{j=1;j\neq r}^N P(R_r=m,R_j=n|U_r=u)&=&(N-1)h_{mn}(u)
=P(R_r=m|U_r=u)\\
&=&{N-1\choose m-1}u^{m-1}(1-u)^{N-m}=(N-1)h_m(u)
\end{eqnarray*}
\begin{eqnarray*}
\sum_{i=1;i\neq r}^N P(R_i=m,R_r=n|U_r=u)&=&(N-1)h_{nm}(u)
=P(R_r=n|U_r=u))\\
&=&{N-1\choose n-1}u^{n-1}(1-u)^{N-n}=(N-1)h_n(u)
\end{eqnarray*}
Thus for $0<u<1$ and $m\neq n$
\[
\sum_{i\neq j=1}^N P(R_i=m,R_j=n|U_r=u)=
(N-1)(N-2)g_{mn}(u)+(N-1)[h_m(u)+h_n(u)]
=1
\]
and hence
\[
g_{mn}(u)=\frac{1}{(N-1)(N-2)}\left[
1-(N-1)[h_{m}(u)+h_{n}(u)]\rule{0in}{.2in}
\right]
\]
and
\[
h_{m}(u)=\frac{1}{N-1}{N-1\choose m-1}u^{m-1}(1-u)^{N-m}
\]

\begin{eqnarray*}
E(J_{0i}J_{\ell j}|U_r)&=&
P(R_i\in A_0, R_j\in A_\ell|U_r)
=\sum_{m\in A_0}\sum_{n\in A_\ell}P(R_i=m,R_j=n|U_r)\\
&=& \s\sum_{m\in A_0}\sum_{n\in A_\ell} g_{mn}(U_r)\quad\mbox{when $i,j,r$ are distinct}\\
&=& \s\sum_{m\in A_0}\sum_{n\in A_\ell} h_{m}(U_r)= 
n_\ell\sum_{m\in A_0}h_{m}(U_r)\quad\mbox{when $i=r\neq j$}\\
&=&\s \sum_{m\in A_0}\sum_{n\in A_\ell} h_{n}(U_r)=
n_0\sum_{n\in A_\ell} h_{n}(U_r)\quad\mbox{when $j=r\neq i$}
\end{eqnarray*}
Making use of (\ref{psi0}), (\ref{psi1}) and (\ref{psi2})
we have
\[
\sum_{i\neq j=1; i,j\neq r}^N
\hat{\psi}_{ij}(v_N)=\sum_{i,j}
\hat{\psi}_{ij}(v_N)- \hat{\psi}_{r\cdot}-\hat{\psi}_{\cdot r}=0
\]
and thus
\begin{eqnarray*}
E(T_\ell|U_r)&=&
\sum_{i\neq j=1; i,j\neq r}^N
\hat{\psi}_{ij}(v_N)E(J_{0i}J_{\ell j}|U_r)+
\sum_{j=1; j\neq r}^N
\hat{\psi}_{kj}(v_N)E(J_{0r}J_{\ell j}|U_r)\\
&&\hspace{1in}\qquad
+
\sum_{i=1; i\neq r}^N
\hat{\psi}_{ik}(v_N)E(J_{0i}J_{\ell r}|U_r)\\
&=&\sum_{m\in A_0}\sum_{n\in A_\ell} g_{mn}(U_r)
\sum_{i\neq j=1; i,j\neq r}^N\hat{\psi}_{ij}(v_N)  +n_\ell
\sum_{m\in A_0} h_{m}(U_r)\hat{\psi}_{r\cdot}(v_N)\\
&&\hspace{1in}\qquad +
n_0\sum_{n\in A_\ell} h_{n}(U_r)
\hat{\psi}_{\cdot r}(v_N)\\
&=&
(v_{Nr}-\bar{v}_N)\left[
n_0
\sum_{n\in A_\ell} h_{n}(U_r)-n_\ell
\sum_{m\in A_0} h_{m}(U_r)
\right]=(v_{Nr}-\bar{v}_N) Z_{\ell r}
\end{eqnarray*}
where
\begin{eqnarray*}
Z_{\ell r}&=&\frac{n_0n_\ell}
{N-1}\left[\frac{1}{n_\ell}\sum_{n\in A_\ell} {N-1\choose n-1}U_r^{n-1}(1-U_r)^{N-n}\right.\\
&&\qquad \left.
-\frac{1}{n_0}\sum_{m\in A_0} {N-1\choose m-1}U_r^{m-1}(1-U_r)^{N-m}
\right]\\
&=&\frac{n_0n_\ell}{N-1}
\left[
\frac{1}{n_\ell}E(I_{A_\ell}(X)|U_r)
-\frac{1}{n_0}E(I_{A_0}(X)|U_r)
\right] \\
&=&\frac{n_0n_\ell}{N-1}
E(B_\ell|U_r)
\end{eqnarray*}
with 
\[
P(X=n|U)={N-1\choose n-1}U^{n-1}(1-U)^{N-n}\quad\mbox{\quad for $n=1,\ldots,N$}
\]
and 
\[
B_\ell=\frac{1}{n_\ell}I_{A_\ell}(X)-\frac{1}{n_0}I_{A_0}(X)
\]
The unconditional distribution of $X$ is
\[
P(X=n)=\int_0^1{N-1\choose n-1}u^{n-1}(1-u)^{N-n}du=\frac{1}{N}
\quad\mbox{\quad for $n=1,\ldots,N$}
\]
The $Z_{\ell r}$ are bounded and iid random variables with mean 0 and variance
$\tau^2=\tau^2(n_0,n_\ell,N)$. From the Lindeberg central limit theorem
it follows that 
\[
\frac{\sum_{r=1}^N (v_{Nr}-\bar{v}_N)Z_{\ell r}}{
\sqrt{\sum_{r=1}^N (v_{Nr}-\bar{v}_N)^2
\tau^2(n_0,n_\ell,N)}}\stackrel{{\cal L}}{\longrightarrow} {\cal N}(0,1)
\]
provided that
\begin{equation}
\frac{\max_r(v_{Nr}-\bar{v}_N)^2}{\sum_{r=1}^N(v_{Nr}-\bar{v}_N)^2}
\rightarrow 0\quad
\mbox{as \quad $n_0, n_\ell, N\rightarrow \infty$}
\label{vratio}
\end{equation}
(\ref{vratio}) was shown to hold in Example 17 of page 355 in
\cite{Lehmann06} under $C_2$.

In order to conclude that 
\[
\frac{T_\ell}{\sqrt{\var(T_\ell)}}=\frac{W_{X,Y_\ell}^\ast-n_0n_\ell/2}{
\sqrt{\var(W_{X,Y_\ell}^\ast)}}\stackrel{{\cal L}}{\longrightarrow} {\cal N}(0,1)
\]
it suffices (\cite{Lehmann06}, p.~349 and p.~363)
to show that
\begin{equation}
\frac{\var(W_{X,Y_\ell}^\ast)}{\sum_{r=1}^N (v_{Nr}-\bar{v}_N)^2
\tau^2(n_0,n_\ell,N)}=
\frac{\var(W_{X,Y_\ell}^\ast)}{N\sigma^2(v_N)
\tau^2(n_0,n_\ell,N)}
\longrightarrow 1
\label{ratio}
\end{equation}
as\quad $n_0, n_\ell, N\rightarrow \infty$

From \cite{Lehmann06}, p.~355, we have
\[
\sigma^2(v_N)=\frac{1}{N}\sum_{i=1}^N(v_{Ni}-\bar{v}_N)^2=
\frac{N^2-1}{12}-
\frac{\sum_{i=1}^e(d_i^3-d_i)}{12N}
\]
or
\[
\frac{\sum_{i=1}^e(d_i^3-d_i)}{N}=N^2-1-12\sigma^2(v_N)
\]
and thus using (\ref{varformula})
\[
\var(W_{X,Y_\ell}^\ast)=
\frac{n_0 n_\ell}{12}\left[
n_0+n_\ell+1
-\frac{n_0+n_\ell-2}{(N-1)(N-2)}\left(
N^2-1-12\sigma^2(v_N)
\right)\right.
\]
\[
\hspace{3in}\left.
-\frac{3(N-n_0-n_\ell)}{N(N-1)(N-2)}\sum_{i=1}^e d_i(d_i-1)
\right]
\]
With
\[
n_0+n_\ell+1
-\frac{n_0+n_\ell-2}{(N-1)(N-2)}\left(
N^2-1\right)=\frac{3(N-n_0-n_\ell)}{N-2}
\]
we have
\[
\var(W_{X,Y_\ell}^\ast)=
\frac{n_0 n_\ell}{12}\left[
\frac{3(N-n_0-n_\ell)}{N-2}+\frac{n_0+n_\ell-2}{(N-1)(N-2)}12\sigma^2(v_N)
\right.
\]
\[
\left.\hspace{2.5in}
-\frac{3(N-n_0-n_\ell)}{N(N-1)(N-2)}\sum_{i=1}^e d_i(d_i-1)
\right]
\]
Assuming condition $C_2$
we have (see \cite{Lehmann06}, p.~356)
$\sum d_i(d_i-1)\leq N(1-\epsilon)N-N=O(N^2)$ and 
$1/12\geq \sigma^2(v_N)/N^2\geq \epsilon/12$, i.e., 
\begin{equation}
\frac{\var(W_{X,Y_\ell}^\ast)}{n_0n_\ell(n_0+n_\ell)\sigma^2(v_N)/N^2}\rightarrow 1
\quad
\mbox{as \quad $N, n_0, n_\ell\rightarrow\infty$}
\label{varW}
\end{equation}
Thus condition (\ref{ratio}) is equivalent to
\begin{equation}
\frac{n_0n_\ell(n_0+n_\ell)\sigma^2(v_N)/N^2}{N
\sigma^2(v_N)\tau^2(n_0,n_\ell,N)}=
\frac{n_0n_\ell(n_0+n_\ell)}{N^3
\tau^2(n_0,n_\ell,N)}\longrightarrow 1
\label{ratio2}
\end{equation}
as\quad $N, n_0, n_\ell\rightarrow\infty$.
Using the identity
\begin{equation}
\var(B_\ell)=\var(E(B_\ell|U))+E(\var(B_\ell|U))
\label{varidentity}
\end{equation}
we find
\[
\var(B_\ell)=E(B_\ell^2)=\frac{1}{n_{\ell}^2}\frac{n_\ell}{N}+
\frac{1}{n_0^2}\frac{n_0}{N}=\frac{n_0+n_\ell}{Nn_\ell n_0}
\]
and
\[
\tau^2(n_0,n_\ell,N)=\var\left(
\frac{n_0n_\ell}{N-1}E(B_\ell|U_r)
\right)
=\frac{n_0n_\ell (n_0+n_\ell)}{N(N-1)^2}
-\frac{n_0^2n_\ell^2}{(N-1)^2}E(\var(B_\ell|U_r))
\]
Condition (\ref{ratio2}) is equivalent to
\begin{equation}
\frac{N(N-1)^2}{n_0n_\ell(n_0+n_\ell)}
\frac{n_0^2n_\ell^2}{(N-1)^2}E(\var(B_\ell|U_r))=
\frac{n_0n_\ell N}{n_0+n_\ell}E(\var(B_\ell|U_r)) \longrightarrow 0 
\label{ratio3}
\end{equation}
as $N, n_0, n_\ell\rightarrow\infty$
\begin{eqnarray*}
E(B_\ell|U_r)&=&\frac{1}{n_\ell}P(X\in A_\ell|U_r)-\frac{1}{n_0}P(X\in A_0|U_r)\\
\var(B_\ell|U_r) &=& \frac{1}{n_\ell^2}P(X\in A_\ell|U_r)(1-P(X\in A_\ell|U_r))\\
&& \qquad +
\frac{1}{n_0^2}P(X\in A_0|U_r)(1-P(X\in A_0|U_r))\\
&&\qquad +
\frac{2}{n_\ell n_0}P(X\in A_0|U_r)P(X\in A_\ell|U_r)
\end{eqnarray*}
With condition $C_1$ it follows from the Berry-Essen theorem \cite{vBeek72}
that all three terms
in 
$\var(B_\ell|U_r)n_0n_\ell N/(n_0+n_\ell)$
converge to zero at the rate of $1/\sqrt{N}$ uniformly in $U_r$, except for $U_r$ in arbitrarily 
small intervals near 
zero and near $n_0/N$, $m_{\ell-1}/N$,
$m_{\ell}/N$ and $1$.

Note that condition $C_1$ guarantees the boundedness of  
\[
\frac{n_0n_\ell N}{n_0+n_\ell}\left(\frac{1}{n_\ell^2}+
\frac{1}{n_0^2}+\frac{2}{n_0n_\ell}\right)=\frac{n_0n_\ell N}{n_0+n_\ell}\left(\frac{n_0+n_\ell}{n_0n_\ell}\right)^2=\frac{N(n_0+n_\ell)}{n_0n_\ell}
\]
Thus the portion of $E(\var(B_\ell|U_r))n_0n_\ell N/(n_0+n_\ell)$ over the small
intervals can be held below any $\epsilon>0$, simply by the smallness
of those intervals, since $N(n_0+n_\ell)/(n_0n_\ell)$ stays bounded.
The remainder portion of $E(\var(B_\ell|U_r))n_0n_\ell N/(n_0+n_\ell)$ is of order $1/\sqrt{N}$. 
Thus condition (\ref{ratio3}) is satisfied
as $N, n_0, n_\ell\rightarrow\infty$.

To show asymptotic multivariate normality of
$(W_{X,Y_1}^\ast, \ldots, W_{X,Y_K}^\ast)$ it suffices 
(according to the Cram\'er-Wold device) to show
asymptotic univariate normality of
\[
S=\sum_{\ell=1}^K \lambda_\ell (W_{X,Y_\ell}^\ast-n_0n_\ell/2)=
\sum_{\ell=1}^K\lambda_\ell T_\ell\quad
\mbox{for any $\lambda_1,\ldots, \lambda_K$}
\]
We may assume that $\sum_{i=1}^K\lambda_i^2\neq 0$, the 
$=0$ case being trivial.
Using again the projection approach on $S$, i.e., approximate
$S$ by
\[
\hat{S}=\sum_{r=1}^N E(S|U_r)=\sum_{r=1}^N  \sum_{\ell=1}^K \lambda_\ell E(T_\ell|U_r)
=\sum_{r=1}^N(v_{Ni}-\bar{v}_N)
\sum_{\ell=1}^K\lambda_\ell \frac{n_0 n_\ell}{N-1}E(B_\ell|U_r)
\]
with
\[
B_\ell=\frac{1}{n_\ell}I_{A_\ell}(X)-\frac{1}{n_0}I_{A_0}(X)\quad
\mbox{and $X$ as before.}
\]
The asymptotic normality of the $\hat{S}/\sqrt{\var(\hat{S})}$ follows again from 
the Lindeberg central limit theorem provided that (\ref{vratio}) holds,
which is a consequence of condition $C_2$. The  asymptotic normality of the $S/\sqrt{\var(S)}$ then follows if we are able to show
$\var(S)/\var(\hat{S})\rightarrow 1$ as 
$n_0,\ldots, n_K\rightarrow \infty$.

As in (\ref{varW}) we have for $\ell\neq k$
\[
\frac{\cov(T_\ell,T_k)}{n_0 n_\ell n_k\sigma^2(v_N)/N^2}
\rightarrow 1 \quad
\mbox{or}\quad
\frac{\cov(T_\ell,T_k)}{N\sigma^2(v_N)}\rightarrow
\gamma_0\gamma_\ell \gamma_k\quad
\mbox{and}\quad
\frac{\var(T_\ell)}{N\sigma^2(v_N)}\rightarrow
\gamma_0\gamma_\ell (\gamma_0+\gamma_\ell)
\]
as $n_0,n_1, \ldots, n_K\rightarrow \infty$.
Thus
\begin{eqnarray*}
\frac{\var(S)}{N\sigma^2(v_N)}&\rightarrow&
{\sum\sum}_{\ell\neq k}\lambda_\ell\lambda_k\gamma_0\gamma_\ell\gamma_k
+\sum_{\ell=1}^K\lambda_\ell^2 \gamma_0\gamma_\ell(\gamma_0+\gamma_\ell)\\
&& = \gamma_0(1-\gamma_0)^2
\left(\sum_{\ell=1}^K \lambda_\ell \gamma_\ell^\ast \right)^2
+\gamma_0(1-\gamma_0)\left(
\var(\lambda)+\left(\sum_{\ell=1}^K \lambda_\ell \gamma_\ell^\ast \right)^2
\right)=\omega^2 
\end{eqnarray*}
where the $\gamma_\ell^\ast=\gamma_\ell/(1-\gamma_0), \ell=1,\ldots,K$,
sum to one. From our assumption
that the $\lambda_\ell$ are not all zero it follows that $\omega^2>0$.

As in (\ref{varidentity}) we have
\begin{eqnarray*}
\cov(B_\ell,B_k)=E(\cov(B_\ell,B_k|U)+\cov(E(B_\ell|U),E(B_k|U))
\end{eqnarray*}
For $\ell\neq k$ we have
$\cov(B_\ell,B_k)=1/(n_0 N)$ and 
\begin{eqnarray*}
\cov(B_\ell,B_k|U)&=& \frac{1}{n_0^2}P(X\in A_0|U)(1-P(X\in A_0|U))\\
&& -\frac{1}{n_\ell n_k}P(X\in A_\ell|U)P(X\in A_k|U)\\
&& + \frac{1}{n_\ell n_0}P(X\in A_\ell|U)P(X\in A_0|U)\\
&& +\frac{1}{n_0 n_k}P(X\in A_0|U)P(X\in A_k|U)
\end{eqnarray*}
Again it can be shown that  
$E(\cov(B_\ell,B_k|U))n_0^2n_\ell n_k/N^2$
converges to zero 
just as in the case of $E(\var(B_\ell|U))n_0n_\ell N/(n_0+n_\ell)$ previously.

\begin{eqnarray}
\var(\hat{S})&=&\sum_{r=1}^N (v_{rN}-\bar{v}_N)^2\sum_{\ell=1}^K
\sum_{k=1}^K 
\lambda_\ell \lambda_k \frac{n_0^2 n_\ell n_k}{(N-1)^2}
\cov(E(B_\ell|U),E(B_k|U))\nonumber\\
\frac{\var(\hat{S})}{N\sigma^2(v_N)}&=&
\sum_{\ell=1}^K\sum_{k=1}^K 
\lambda_\ell \lambda_k \frac{n_0^2 n_\ell n_k}{(N-1)^2}
\cov(E(B_\ell|U),E(B_k|U))\nonumber\\
&=&
{\sum\sum}_{k\neq \ell}\lambda_\ell \lambda_k \frac{n_0^2 n_\ell n_k}{(N-1)^2}\left(\frac{1}{n_0N}-E(\cov(B_\ell,B_k|U)\right)\nonumber\\
&& \qquad +\sum_\ell \lambda_\ell^2\frac{n_0^2n_\ell^2}{(N-1)^2}\left(
\frac{n_0+n_\ell}{Nn_0n_\ell}-E(\var(B_\ell|U))
\right)\label{last}
\end{eqnarray}
Here
\[
{\sum\sum}_{k\neq \ell}\lambda_\ell \lambda_k \frac{n_0^2 n_\ell n_k}{(N-1)^2}\frac{1}{n_0N} +\sum_\ell \lambda_\ell^2\frac{n_0^2n_\ell^2}{(N-1)^2}
\frac{n_0+n_\ell}{Nn_0n_\ell}
\rightarrow \omega^2
\]
as $n_0, n_1, \ldots, n_K\rightarrow \infty$. The remaining terms 
in (\ref{last}) converge to zero. Thus $\var(S)/\var(\hat{S})
\rightarrow 1$.

\bibliography{Steel}
\end{document}